\documentclass[12pt,oneside,reqno]{amsart}
\usepackage{bbm}
\usepackage{mathrsfs}
\usepackage{stmaryrd}
\usepackage{amsfonts}
\usepackage{enumerate,amsmath,amssymb,amsthm}
\newcommand{\be}{\begin{eqnarray}}
\newcommand{\ee}{\end{eqnarray}}
\newcommand{\ce}{\begin{eqnarray*}}
\newcommand{\de}{\end{eqnarray*}}
\newtheorem{theorem}{Theorem}[section]
\newtheorem{lemma}[theorem]{Lemma}
\newtheorem{remark}[theorem]{Remark}
\newtheorem{definition}[theorem]{Definition}
\newtheorem{proposition}[theorem]{Proposition}

\newtheorem{corollary}[theorem]{Corollary}

\def\a{\alpha}

\def\Om{\Omega}
\def\b{\beta}
\def\p{\partial}

\def\l{\lambda}

\def\[{{\Big[}}
\def\]{{\Big]}}
\def\<{{\langle}}
\def\>{{\rangle}}
\def\({{\Big(}}
\def\){{\Big)}}

\def\dif{{\mathord{{\rm d}}}}

\def\div{{\mathord{{\rm div}}}}

\def\min{{\mathord{{\rm min}}}}

\def\u{\mathord{{\bf u}}}
\def\f{\mathord{{\bf f}}}
\def\v{\mathord{{\bf v}}}

\def\w{\mathord{{\bf w}}}
\def\e{\mathord{{\bf e}}}

\def\b{\mathord{{\bf b}}}

\def\no{\nonumber}
\def\bt{\begin{theorem}}
\def\et{\end{theorem}}
\def\bl{\begin{lemma}}
\def\el{\end{lemma}}
\def\br{\begin{remark}}
\def\er{\end{remark}}
\def\bd{\begin{definition}}
\def\ed{\end{definition}}
\def\bp{\begin{proposition}}
\def\ep{\end{proposition}}
\def\bc{\begin{corollary}}
\def\ec{\end{corollary}}
\def\bC{{\mathbf C}}
\def\cA{{\mathcal A}}
\def\cB{{\mathcal B}}

\def\cO{{\mathcal O}}

\def\cQ{{\mathcal Q}}

\def\cU{{\mathcal U}}

\def\mN{{\mathbb N}}

\def\mR{{\mathbb R}}

\def\mT{{\mathbb T}}

\def\mX{{\mathbb X}}

\def\sD{{\mathscr D}}
\def\sE{{\mathscr E}}

\def\sP{{\mathscr P}}

\def\geq{\geqslant}
\def\leq{\leqslant}

\def\bH{{\mathbf H}}

\def\bL{{\mathbf L}}
\def\bW{{\mathbf W}}

\def\Om{{\Omega}}
\def\lb{{\llbracket}}
\def\rb{{\rrbracket}}

\def\U{{\mathbf U}}
\allowdisplaybreaks

\begin{document}

\title{A Tamed 3D Navier-Stokes Equation in Domains}

\date{}
\author{Xicheng Zhang }


%Galerkin's Approximation, Analytic Semigroup,

\dedicatory{
School of Mathematics and Statistics\\
The University of New South Wales, Sydney, 2052, Australia\\
Department of Mathematics,
Huazhong University of Science and Technology,\\
Wuhan, Hubei 430074, P.R.China\\
Email: XichengZhang@gmail.com
 }

\begin{abstract}

In this paper, we analyze a tamed  3D Navier-Stokes equation in uniform $C^2$-domains
(not necessarily bounded), which obeys the scaling invariance principle,
and prove the existence and uniqueness of strong solutions to this tamed equation.
In particular, if there exists a bounded solution to the classical
3D Navier-Stokes equation, then this solution satisfies our tamed equation.
Moreover, the existence of  a global attractor for the tamed equation in bounded domains
is also proved. As simple applications, some well known results for the classical Navier-Stokes
equations in unbounded domains are covered.
\end{abstract}

\maketitle \rm

\section{Introduction}

The motion of a viscous incompressible fluid in a domain $\Om\subset\mR^3$
is described by the Navier-Stokes equation (NSE) as follows (with homogeneous boundary):
\be\label{Ns1}
\left\{
\begin{aligned}
&\p_t\u =\nu\Delta \u -(\u \cdot\nabla)\u +\nabla P+\f,\\
&\div(\u)=0,
\ \ (t,x)\in[0,\infty)\times\Om,\\
&\u(t,x)=0,\ \ t\geq 0,\ \ x\in\p\Om,\ \ \u(0)=\u_0,
\end{aligned}
\right.
\ee
where $\nu>0$ is the kinematic viscosity constant,
$\u(t,x)=(u_1(t,x),u_2(t,x),u_3(t,x))$ represents the velocity field,
$P=P(t,x)$ is the pressure (an unknown scalar function),
$\f$ is a known external force.

The study of 3D NSEs has a long history. In their pioneering works,
Leray \cite{Le} and Hopf \cite{Ho} proved the existence of a weak solution
to equation (\ref{Ns1}). Since then, there are many papers devoted to the study of regularities of
 Leray-Hopf weak solutions (cf. \cite[etc.]{La, Te, So}).
 Up to now, one knows that the singular set of the Leray-Hopf
weak solutions has Lebesgue measure zero (cf. \cite{Le, He, Ga}). Moreover,
a deep result obtained by Scheffer \cite{Sc}
and Caffarelli, Kohn and Nirenberg \cite{Ca-Ko-Ni} says that the singular set for a class of
weak solutions (satisfying a generalized energy inequality) has one
dimensional Hausdorff measure zero (see also \cite{Li}).
However, the uniqueness and regularity of Leray-Hopf weak solutions
are still big open problems.

Most of the source of difficulties to solve equation (\ref{Ns1})
comes from the nonlinear term $(\u\cdot\nabla)\u$ (cf. \cite{Ga}).
In order to counteract this term, the authors in \cite{Ro-Zh}
analyzed the following modified (called tamed therein)
3D NSE in   $\Om=\mR^3$:
\be\label{Ns2}
\left\{
\begin{aligned}
&\p_t \u =\nu\Delta \u -(\u \cdot\nabla)\u +\nabla P -g^\nu_N(|\u |^2)\u+\f,\\
&\div(\u)=0,\ \ t\geq 0,\ \ x\in\mR^3,\ \ \u(0)=\u_0,
\end{aligned}
\right.
\ee
where $|\u|^2:=\sum^3_{j=1}|u_j|^2$ and for $N>0$
\be
g^\nu_N(r):=(r-N)\cdot 1_{\{r\geq N\}}/\nu.\label{Con}
\ee
The existence of a unique smooth solution to equation (\ref{Ns2}) was proved in \cite{Ro-Zh}
when the initial velocity is smooth (in Sobolev spaces).
The main feature of equation (\ref{Ns2}) is that if there exists a bounded solution (say bounded
by $\sqrt{N}$ for some large $N$) to the classical NSE, then this solution
must satisfy equation (\ref{Ns2}). Therein, the property
that the Leray projection operator onto the space of divergence free vector fields
commutes with the derivatives plays a key role. But,
when we consider NSE (\ref{Ns1}) in a domain, this property does not hold in general (cf. \cite[p.83-85]{Lions}).

In order to deal with the Dirichlet boundary problem and keep the same feature as equation (\ref{Ns2}),
in the present paper, we consider the following globally tamed scheme (assuming $\f=0$ for simplicity):
\be
\label{Tamed}
\p_t \u=\nu\Delta \u-(\u\cdot\nabla)\u +\nabla P-g^{\nu,\kappa}_N(\|\u-\U\|^2_\infty)(\u-\U),
\ee
where $\|\u\|_\infty:=\sup_{x\in\Om}|\u(x)|$, $\U$ is a reference  velocity field
and for $\kappa,N\geq 1$
$$
g^{\nu,\kappa}_{N}(r):=\kappa\cdot(r-N)1_{\{r\geq N\}}/\nu.
$$
Here, $\kappa\geq 1$ is a dimensionless constant and $\sqrt{N}$
has the velocity dimension.

Let $(\u_{N,\U},P_{N,\U})$ be a solution pair of equation (\ref{Tamed}). Simple calculations show that
$(\u_{N,\U},P_{N,\U})$ has the following properties:
\begin{enumerate}[{\bf (A)}]
\item (Galilean invariance): for any constant velocity vector $\v\in\mR^3$
\ce
\u^{\v}_{N,\U}(t,x)&:=&\u_{N,\U+\v}(t,x-\v t)+\v,\\
P^{\v}_{N,\U}(t,x)&:=&P_{N,\U+\v}(t, x-\v t)
\de
is also a solution pair of equation (\ref{Tamed}).
\item (Rotation symmetry): for any orthogonal matrix $\cQ$ (i.e. $\cQ\cQ^t=I$)
\ce
\u^{\cQ}_{N,\U}(t,x)&:=&\cQ^t\u_{N,\cQ^t\U}(t,\cQ x),\\
P^\cQ_{N,\U}(t,x)&:=&P_{N,\cQ^t\U}(t,\cQ x)
\de
is also a solution pair of equation (\ref{Tamed}).
\item (Scale invariance): for any $\lambda>0$
\ce
\u^{\lambda}_{N,\U}(t,x)&:=&\lambda\u_{\lambda^{-2}N,\lambda^{-1}\U}(\lambda^2t,\lambda x),\\
P^\lambda_{N,\U}(t,x)&:=&\lambda^2P_{\lambda^{-2}N,\lambda^{-1}\U}(\lambda^2t,\lambda x)
\de
is also a solution pair of equation (\ref{Tamed}).
\end{enumerate}
These three properties are exhibited by the classical Navier-Stokes equations  (cf. \cite{Be-Ma}).
%Although the term $-g^{\nu,\kappa}_N(\|\u\|^2_\infty)\u$ is artificial, equation (\ref{Tamed})
%has many similarities with the classical NSEs. Moreover, as we can see below,
%all the difficulties met in equation (\ref{Ns1}) does not appear for equation (\ref{Tamed}).

Intuitively, when the maximum of the fluid velocity is larger than $\sqrt{N}$, the dissipative term
$g^{\nu,\kappa}_N(\|\u\|^2_\infty)\u$ (regarded as some extra force) will enter into the equation and
restrain the flux of the liquid. In this sense, the value of $N$ plays the role of a valve.
On the other hand, when we realize equation (\ref{Tamed}) on a computer,
the value of $N$ can be reset as an arbitrarily large number along with the process of calculations
as long as there is no explosion. So, the term involving
$g^{\nu,\kappa}_N$ plays the role of some kind of adjustment.
The parameter $\kappa$ can be understood as the extent of the extra dissipative force,
and will be used to give a better estimate for $\|\u\|_\infty$ in terms of $N$
(see part (III) of Theorem \ref{main}).

In contrast with equation (\ref{Ns2}), the tamed equation (\ref{Tamed}) in domain $\Om$
is global since $g^{\nu,\kappa}_N(\|\u\|^2_\infty)$ depends on all values of $\u$ in $\Om$.
But, better than (\ref{Ns2}), it is easy to write down the vorticity equation:
Let $\omega=\mathrm{curl}\u=\nabla\wedge\u$. Then
$$
\p_t\omega=\nu\Delta\omega+(\omega\cdot\nabla)\u-(\u\cdot\nabla)\omega
-g^{\nu,\kappa}_N(\|\u-\U\|^2_\infty)\omega.
$$
We remark that in \cite{Ca-Re-Kl}, Caraballo, Real and Kloeden
studied the following globally modified NSE in a bounded regular domain $\Omega$:
\be\label{Ns0}
\left\{
\begin{aligned}
&\p_t \u(t)=\Delta \u-\min\{1,N/\|\nabla\u\|_{\bL^2}\}(\u\cdot\nabla)\u+\nabla P,\\
&\div(\u)=0,
\ \ (t,x)\in[0,\infty)\times\Om,\\
&\u(t,x)=0,\ \ t\geq 0,\ \ x\in\p\Om,\ \ \u(0)=\u_0,
\end{aligned}
\right.
\ee
and they proved the existence of a unique strong solution to
this modified equation as well as the existence of a global attractor.
Nevertheless, equation (\ref{Ns0}) does not enjoy the above properties
{\bf (A)}-{\bf (C)}.

This paper is organized as follows: In Section 2, all main results are announced.
In Section 3, we prepare some necessary lemmas for later use. In the remaining sections,
we shall give the proofs of main results. We want to emphasize that
for the proof of existence of strong solutions (see Section 4),
not using the usual Galerkin approximation, we only use the linearized equations and
simple Picard's iteration. Moreover, the semigroup method used in Fujita-Kato \cite{Fu-Ka}
(cf. \cite{So}) will be used to improve the regularity of strong solutions (see Section 5).
The existence of a global attractor for the evolution semigroup determined by equation (\ref{Tamed})
will follow by proving some asymptotic compactness (cf. \cite[etc.]{Ca-Re-Kl, Te2}).

\section{Announcement of Main Results}

Throughout this paper, all $\mR^3$-valued functions and spaces of such functions will
be denoted by boldfaced letters, and we use the following convention:
the letter $C$ with or without subscripts will
denote a positive constant whose value may change in different occasions.

Let $\Om$ be a uniform $C^3$-regular domain of $\mR^3$ (see \cite[p.84]{Ad} for the definition
of regular domains).
Let $\bC^\infty_0(\Om)$ denote the set of all smooth functions from $\Om$ to $\mR^3$ with
compact supports in $\Om$, and $\bC^\infty_{0,\sigma}(\Om)\subset \bC^\infty_0(\Om)$
the set of all smooth vector fields of divergence free. For $p>1$, let
$\bL^p(\Om)$ be the usual $\mR^3$-valued $L^p$-space with the norm denoted by
$\|\cdot\|_{\bL^p(\Omega)}=\|\cdot\|_{\bL^p}$,
and $\bL^2_\sigma(\Om)$  the closure of $\bC^\infty_{0,\sigma}(\Om)$ in $\bL^p(\Om)$.
For $k\in\mN$ and $p>1$, let $\bW^{k,p}(\Om)$ be the space of $\mR^3$-valued functions with finite norm:
\ce
\|\u\|_{\bW^{k,p}(\Om)}:=\|\u\|_{\bW^{k,p}}:=\bigg(\sum_{j=0}^k\int_\Om|\nabla^j\u(x)|^p\dif x\bigg)^{\frac{1}{p}}<+\infty,
\de
where $\nabla^j$ denotes the $j$-th order generalized derivative operator. The space
$\bW^{1,2}_0(\Om)$ (resp. $\bW^{1,2}_{0,\sigma}(\Om)$) denotes the completion
of $\bC^\infty_0(\Om)$ (resp. $\bC^\infty_{0,\sigma}(\Om)$)
with respect to the above norm with $k=1$ and $q=2$.

Let $\sP$ be the orthogonal projection  from $\bL^2(\Om)$ to $\bL^2_\sigma(\Om)$.
By $A$ (called the Stokes operator) we denote the self-adjoint operator in $\bL^2_\sigma(\Omega)$
formally given by
$$
A:=-\sP\Delta.
$$
More precisely, $\u\in\sD(A)$ if and only if
for some $\w\in\bL^2_\sigma(\Om)$ (written as $A\u=\w$), it holds that
$$
\<\nabla \u,\nabla\v\>_{\bL^2}=\<\w,\v\>_{\bL^2},\quad\forall \v\in\bW^{1,2}_{0,\sigma}(\Om).
$$
In particular, $\sD(A^{\frac{1}{2}})=\bW^{1,2}_{0,\sigma}(\Om)$ and
\be
\|A^{\frac{1}{2}}\u\|_{\bL^2}=\|\nabla\u\|_{\bL^2},\quad \u\in\bW^{1,2}_{0,\sigma}(\Om).
\label{Id}
\ee
Moreover, it is well known that (cf. \cite[p.129]{So})
$$
\sD(A)=\bW^{2,2}(\Om)\cap\bW^{1,2}_{0,\sigma}(\Om).
$$
Since $A$ is a positive self-adjoint operator in $\bL^2_\sigma(\Om)$,
for $\a\in(-1,1)$, the fractional power $A^\a$ is well defined via the spectral representation.
For $\beta\in[0,2]$, define the Hilbert space
$$
\bH^{\beta}(\Omega):=\bH^{\beta}:=\sD(A^{\beta/2})
$$
with the norm $\|\cdot\|_{\bH^\beta}$ generated by inner product
$$
\<\u,\v\>_{\bH^\beta}:=\<\u,\v\>_{\bL^2}+\<A^{\beta/2}\u,A^{\beta/2}\v\>_{\bL^2}.
$$

We introduce the following bilinear form $B$ on $\bW^{1,2}_{0,\sigma}(\Om)=\bH^1$:
\be
B(\v,\u):=-\sP((\v\cdot\nabla)\u).\label{defB}
\ee
Using $\sP$ to act on both sides of equation (\ref{Tamed}),
we can and shall consider the following equivalent abstract equation
\be
\label{ren}
\p_t \u=-\nu A\u+B(\u,\u)-g^{\nu,\kappa}_N(\|\u\|^2_{\infty})\u,\ \ \u(0)=\u_0.
\ee

We give the following definition of strong solutions to the above equation.
\bd\label{def}
Let $\u_0\in\bH^1$. A continuous function
$$
\mR_+\ni t\mapsto\u(t)\in\bH^1
$$
is called a strong solution of equation (\ref{ren}) if
$\u\in L^2_{loc}(\mR_+;\bH^2)$ and for all $t\geq 0$
\be
\u(t)=\u_0-\nu\int^t_0A\u\dif s+\int^t_0B(\u,\u)\dif s
-\int^t_0g^{\nu,\kappa}_N(\|\u\|^2_{\infty})\u\dif s\ \ \mbox{ in \ $\bL^2_\sigma(\Omega)$.}\label{Op1}
\ee
\ed
Our first main result is stated as follows:
\bt\label{Main}
Let $\Om$ be a uniform  $C^2$-domain of $\mR^3$.
For any $\u_0\in\bH^1$, there exists a unique
strong solution $\u(t)=\u_N(t)$ to equation (\ref{ren})
in the sense of Definition \ref{def}, which satisfies
that for any $t\geq 0$
\be
\|\u(t)\|^2_{\bL^2}+2\nu\int^t_0\|\nabla\u\|^2_{\bL^2}
\dif t&\leq&\|\u_0\|^2_{\bL^2},\label{Es22}\\
\|\nabla\u(t)\|^2_{\bL^2}+\nu\int^t_0\|A\u\|^2_{\bL^2}\dif t
&\leq& \frac{\kappa N}{\nu^2}\|\u_0\|^2_{\bL^2}+\|\nabla\u_0\|^2_{\bL^2}\label{Es222}
\ee
and for some $T^*=T^*(\nu,\Omega,\|\u_0\|_{\bL^2})$ and $C=C(\nu,\Omega,\|\u_0\|_{\bL^2})$
\be
\|\nabla\u(t)\|_{\bL^2} \leq C/\sqrt{t},\ \ \forall t\geq T^*.\label{decay}
\ee
Moreover, letting $\u_N(t)$ (resp. $\v_M(t)$) be the solution of equation (\ref{ren}) with initial value
$\u_0\in\bH^1$ (resp. $\v_0\in\bH^1$) and taming function $g^{\nu,\kappa}_N$
(resp. $g^{\nu,\kappa}_M$), we have for any $T>0$
\be
&&\sup_{t\in[0,T]}\|\u_N(t)-\v_M(t)\|_{\bH^1}^2+\int^T_0\|\u_N-\v_M\|_{\bH^1}^2\dif s\no\\
&&\quad\leq C({\nu,N,M,\|\u_0\|_{\bH^1},\|\v_0\|_{\bH^1},T})
\cdot(|N-M|^2+\|\u_0-\v_0\|^2_{\bH^1}),\label{Col}
\ee
where the constant $C({\nu,N,M,\|\u_0\|_{\bH^1},\|\v_0\|_{\bH^1},T})$
continuously depends on its parameters.
\et
\br
\rm For $T>0$ and $N\geq 1$, define
$$
\mT^T_N:=\{t\in[0,T]: \|\u(t)\|_{\infty}\geq\sqrt{N}\}.
$$
By  (\ref{Es22}), (\ref{Es222}) and (\ref{PI2}) below, we have
\ce
\l(\mT^T_N)&\leq& \frac{1}{N}\int^T_0\|\u\|^2_\infty\dif s
\leq \frac{C_\Om }{N}\int^T_0\|\u\|_{\bH^2}\cdot\|\nabla\u\|_{\bL^2}\dif s\leq\\
&\leq& \frac{C_\Om}{N}\left(\int^T_0\|\u\|^2_{\bH^2}\dif s\right)^{1/2}
\cdot\left(\int^T_0\|\nabla\u\|_{\bL^2}^2\dif s\right)^{1/2}\\
&\leq& \frac{C_\Om\cdot\|\u_0\|_{\bL^2}\cdot
(\frac{\kappa N}{\nu^2}\|\u_0\|^2_{\bL^2}+\|\nabla\u_0\|^2_{\bL^2}
+\nu T\|\u_0\|^2_{\bL^2})^{1/2}}{\sqrt{2}\nu N},
\de
where $\l(\mT^T_N)$ denotes the Lebesgue measure of $\mT^T_N$.  This gives an estimate
of the length of the time for which $\u_N$ does not satisfy equation (\ref{Ns1}). In particular,
$$
\lim_{N\to\infty}\l(\mT^T_N)=0,
$$
which shows that as $N$ goes to infinity, $\u_N$ satisfies equation (\ref{Ns1}) at ``almost all''
times.
\er
We are now interested in the estimation of $\|\u\|_{\infty}$ in terms of $N$
and prove the following result.
\bt\label{main}
Let $\Om$ be a uniform  $C^2$-domain and $\u_0\in\bH^2$. Let $\u^\kappa_N$
be the unique strong solution in Theorem \ref{Main}. We have the following conclusions:
\begin{enumerate}[{\bf (I)}]
\item There exist two  continuous function $K_1:\mR_+^2\to\mR_+$
and $K_2:\mR_+^3\to\mR_+$ such that for all $t\geq 0$ and $N\geq 1$
\be
\|\u^\kappa_N(t)\|_{\infty}\leq K_1(t,\|\u_0\|_{\bH^2})
+K_2(t,\nu,\|\u_0\|_{\bL^2})\cdot N^3,\label{Es1}
\ee
where $K_1(t,r), K_2(t,\nu,r)\to 0$ as $t\to 0$ or $\nu\to\infty$ or $r\to 0$.
In particular, for $T>0$, if one of the following conditions is satisfied, then there is a
unique strong solution in $[0,T]$ for equation (\ref{Ns1}):
$$
\mbox{(i) $T$ is small;\ \ (ii) $\|\u_0\|_{\bH^2}$ is small;\ \ (iii) $\nu$ is large.}
$$

\item Let $\Om=\mR^3$ or be a bounded uniform $C^4$-domain and $\u=\u^\kappa_N$. Then
$$
\u\in C([0,\infty)\times\bar\Om;\mR^3)
$$
and for $i,j=1,2,3$
$$
\p_t\u,~\p_i\u,~\p_i\p_j\u \in C((0,\infty)\times\bar\Om;\mR^3).
$$
Moreover, for some $P\in C((0,\infty)\times\bar\Om;\mR)$ (with $\int_\Om P(x)\dif x=0$),
it holds that
\be
\p_t \u=\Delta \u
-(\u\cdot\nabla)\u+\nabla P-g_N(\|\u\|^2_\infty)\u,\ \forall (t,x)\in(0,\infty)\times\Om.\label{Cl}
\ee
\item Let $\Omega=\mR^3$ and $\nu>0$. For any $\a>\frac{1}{2}$, there exist $\kappa>0$
and two functions $K_{1,\a,\kappa},K_{2,\a,\kappa}$ as in {\bf (I)}
such that for all $t\geq 0$ and $N\geq 1$
\be
\|\u^\kappa_N(t)\|_{\infty}\leq K_{1,\a,\kappa}(t,\|\u_0\|_{\bH^2})
+K_{2,\a,\kappa}(t,\nu,\|\u_0\|_{\bL^2})\cdot N^\a.\label{Es}
\ee
\end{enumerate}
\et
\br\rm
We do not know whether the $\a$ in (III) can be smaller than $1/2$. If this can be proven,
then (\ref{Ns1}) will have a classical solution. In fact, even for $\a=1/2$,
it seems also hard to prove (\ref{Es}).
\er
\br\rm
Fix $T>0$ and $N_1\geq 1$. Define a sequence of real numbers recursively as follows:
$$
N_{k+1}:=\sup_{t\in[0,T]}\|\u_{N_k}(t)\|^2_\infty,\ \ k\in\mN.
$$
It is easy to see that equation (\ref{Ns1}) has a explosion solution in $[0,T]$
if and only if
$$
N_1<N_2<N_3<\cdots <N_k\rightarrow\infty.
$$
The strict monotonicity is clear. Assume that $\lim_{k\rightarrow\infty}N_k=N_\infty<\infty$.
By the continuous dependence of $\u_N$ with respect to $N$ (see (\ref{Col})), we have
$$
\lim_{k\rightarrow\infty}\sup_{t\in[0,T]}\|\u_{N_k}(t)-\u_{N_\infty}(t)\|^2_\infty=0.
$$
Therefore,
$$
N_\infty=\sup_{t\in[0,T]}\|\u_{N_\infty}(t)\|^2_\infty<\infty,
$$
which implies that $\u_{N_\infty}(t)$ satisfies (\ref{Ns1}), no explosion.
\er

For $\u_0\in\bH^1$, let $\{\u(t;\u_0); t\geq 0\}$ be the unique strong
solution of equation (\ref{ren}), which defines a nonlinear evolution semigroup:
\be
S(t)\u_0:=\u(t;\u_0): \bH^1\to\bH^1.\label{Evo}
\ee
By Theorem \ref{Main}, $\{S(t); t\geq 0\}$ has the following properties:
\begin{enumerate}[(i)]
\item $S(0)=I$ identity map on $\bH^1$;
\item $S(t+s)=S(t)S(s)$ for any $t,s\geq 0$;
\item $[0,\infty)\times\bH^1\ni(t,\u_0)\mapsto S(t)\u_0\in\bH^1$ is continuous.
\end{enumerate}
\bd
A compact subset $\cA\subset\bH^1$ is called a global attractor of
the evolution semigroup $\{S(t); t\geq 0\}$ if
\begin{enumerate}[(i)]
\item $\cA$ is invariant under $S(t)$, i.e, for any $t>0$, $S(t)\cA=\cA$;
\item $\cA$ attracts all bounded set $\cU\subset\bH^1$, i.e.,
$$
\lim_{t\rightarrow\infty}\rho(S(t)\cU,\cA)=0,
$$
where $\rho(\cA_1,\cA_2):=\sup_{\u\in\cA_1}\inf_{\v\in\cA_2}\|\u-\v\|_{\bH^1}$.
\end{enumerate}
\ed

We have the following existence of global attractors of $\{S(t), t\geq 0\}$.
\bt\label{Th3}
Let $\Om$ be a bounded uniform  $C^2$-domain of $\mR^3$.
Then there exists a global attractor
$\cA\subset\bH^1$ to $\{S(t); t\geq 0\}$ defined by (\ref{Evo}).
\et

\section{Preliminaries}

In this section, we collect some necessary materials for later use.
The following lemma is from \cite[Lemma 6]{He}.
\bl\label{Le2}
Let $\phi:\mR_+\to\mR_+$ be an absolute continuous function and
$g:\mR_+\to\mR_+$ a locally Lipschitz continuous function. Suppose that
$\Lambda:=\int^\infty_0\phi(t)\dif t<+\infty$ and $g(\phi)\leq \a\phi^2$ for $\phi\leq\beta$,
where $\a,\beta>0$. If
$$
\phi'(t) \leq g(\phi(t)), \ \forall t\geq 0,
$$
then for $t\geq (\Lambda/\beta)\exp(\a \Lambda)$
$$
\phi(t)\leq (e^{\a  \Lambda}-1)/(\a t).
$$
\el

Let $\{E_\l,\l>0\}$ be the spectrum decomposition of $A$
in $\bL^2_\sigma(\Om)$. The Stokes semigroup is then defined by
$$
e^{-t A}:=\int^\infty_0e^{-t\l}\dif E_\l,
$$
and for $\a\in[-1,1]$, $A^\a$ is given by
$$
A^\a:=\int^\infty_0\l^\a\dif E_\l.
$$

The following lemma is easily derived from the above representations (cf. \cite{So}).
\bl\label{Le8}
(i)
For any $\a\in[0,1]$ and $\u\in\bL^2_\sigma(\Omega)$, we have
$e^{-t A}\u\in \sD(A^\a)$ and
$$
\|A^\a e^{-tA}\u\|_{\bL^2}\leq t^{-\a}\|\u\|_{\bL^2},\ \ \forall t>0.
$$
(ii) For all $\u\in\sD(A^\a)$ and $t\geq 0$
$$
A^\a e^{-t A}\u=e^{-t A}A^\a\u,\ \ \|e^{-tA}\u-\u\|_{\bL^2}\leq C_\a t^\a\|A^\a\u\|_{\bL^2}.
$$
(iii) For any $0\leq\a<\gamma<\beta\leq 1$ and $\u\in \sD(A^\beta)$
$$
\|A^\gamma\u\|_{\bL^2}\leq \|A^\beta\u\|_{\bL^2}^{\frac{\gamma-\a}{\beta-\a}}
\cdot\|A^\a\u\|^{\frac{\beta-\gamma}{\beta-\a}}_{\bL^2}
\leq \frac{\gamma-\a}{\beta-\a}\|A^\beta\u\|_{\bL^2}
+\frac{\beta-\gamma}{\beta-\a}\|A^\a\u\|_{\bL^2}.
$$
\el

We recall the following well known results (cf. \cite[Lemma 2.4.2 (p.142),
Lemma 2.5.2 (p.152) and Lemma 2.4.3 (p.143)]{So}).
\bl
(i) For $\a\in[0,1/2]$ and $q=\frac{6}{3-4\a}$, there exists a constant $C=C(\a,q)>0$ such that
for any $\u\in\bH^{2\a}$
\be
\|\u\|_{\bL^q}\leq C\|A^\a\u\|_{\bL^2}.\label{Int}
\ee
(ii) For $\a\in[0,1/2]$ and  $q=\frac{6}{3+4\a}$,
there exists a constant $C=C(\a,q)>0$ such that
for any $\u\in\bL^q(\Omega)$
\be
\|A^{-\a}\sP \u\|_{\bL^2}\leq C\|\u\|_{\bL^q}.\label{Int1}
\ee
(iii) For $\a\in[1/2,1]$ and $q=\frac{6}{5-4\a}$, there exists a constant $C=C(\Omega,\a,q)>0$ such that
for any $\u\in\bH^{2\a}$
\be
\|\u\|_{\bW^{1,q}}\leq C(\|A^\a\u\|_{\bL^2}+\|\u\|_{\bL^2}).\label{Int2}
\ee
\el
This lemma has the following conclusions.
\bl\label{Le6}
Let $\Om$ be a  uniform $C^2$-domain.
For some $C_\Om>0$ and any $\u\in\bH^2(\Om)$
\be
\|\u\|_{\bL^\infty(\Om)}^2\leq C_\Om\cdot\|\u\|_{\bH^2(\Om)}
\cdot\|\nabla\u\|_{\bL^2(\Om)},\label{PI2}
\ee
and for $\frac{3}{4}<\a\leq 1$, some $C_{\a,\Om}>0$ and any $\u\in\bH^{2\a}(\Om)$
\be
\|\u\|_{\bL^\infty(\Om)}\leq C_{\a,\Om}\cdot(\|A^\a\u\|_{\bL^2(\Om)}+\|\u\|_{\bL^2(\Om)}).\label{PI3}
\ee
\el
\begin{proof}
Since $\Om$ is a uniform $C^2$-domain, by \cite[p. 154, Theorem 5.24]{Ad}
there exists a bounded linear operator $E: \bW^{k,p}(\Om)\mapsto\bW^{k,p}(\mR^3)$
such that $E\u=\u$ a.e. on $\Om$. Recall
the following Gagliardo-Nirenberg inequality (cf. \cite[p.24, Theorem 9.3]{Fr}):
Let $1\leq p,q\leq\infty$ and $\a\in[0,1]$ with $p\not=3$ and
$$
\frac{1}{r}=\a\left(\frac{1}{p}-\frac{1}{3}\right)+(1-\a)\frac{1}{q}.
$$
Then, for some $C=C(r,p,q)$ and all $\u\in\bW^{1,p}(\mR^3)\cap\bL^q(\mR^3)$
\be
\|\u\|_{\bL^r(\mR^3)}\leq C\|\u\|_{\bW^{1,p}(\mR^3)}^\a\cdot\|\u\|_{\bL^q(\mR^3)}^{1-\a}.
\label{GN}
\ee
Thus, by (\ref{Int}) we have
\ce
\|\u\|_{\bL^\infty(\Om)}^2&\leq&\|E\u\|_{\bL^\infty(\mR^3)}^2
\leq C\|E\u\|_{\bW^{2,2}(\mR^3)}\cdot\|E\u\|_{\bL^6(\mR^3)}\\
&\leq&C_\Om\|\u\|_{\bW^{2,2}(\Om)}\cdot\|\u\|_{\bL^6(\Om)}\\
&\leq&C_\Om\|\u\|_{\bH^2(\Om)}\cdot\|A^{\frac{1}{2}}\u\|_{\bL^2(\Om)}
\de
and for $q=\frac{6}{5-4\a}>3$
\ce
\|\u\|_{\bL^\infty(\Om)}&\leq&\|E\u\|_{\bL^\infty(\mR^3)}
\leq C_q\|E\u\|_{\bW^{1,q}(\mR^3)}\leq C_{q,\Om}\|\u\|_{\bW^{1,q}(\Om)}\\
&\stackrel{(\ref{Int2})}{\leq}&C_{\a,\Om}\cdot(\|A^\a\u\|_{\bL^2(\Om)}+\|\u\|_{\bL^2(\Om)}).
\de
The proof is complete.
\end{proof}

The contents below in this section are only used in Section 5.

\bl
For some $C,C_\Om>0$ and all $\u\in\bH^2=\sD(A)$, we have
\be
\|A^{-\frac{1}{4}}B(\u,\u)\|_{\bL^2}&\leq& C\|A^\frac{1}{2}\u\|^2_{\bL^2},\label{PI5}\\
\|B(\u,\u)\|_{\bL^2}&\leq& C_\Om(\|A^\frac{5}{8}\u\|^2_{\bL^2}+\|\u\|^2_{\bL^2}).\label{PI6}
\ee
\el
\begin{proof}
By   H\"older's inequality, we have
\ce
\|A^{-\frac{1}{4}}B(\u,\u)\|_{\bL^2}
\stackrel{(\ref{Int1})}{\leq}C\|(\u\cdot\nabla)\u\|_{\bL^{3/2}}
\leq C_\|\u\|_{\bL^6}\cdot\|\nabla\u\|_{\bL^2}
\stackrel{(\ref{Int})}{\leq} C\|A^\frac{1}{2}\u\|^2_{\bL^2}
\de
and
\ce
\|B(\u,\u)\|_{\bL^2}&\leq& \|(\u\cdot\nabla)\u\|_{\bL^2}
\leq  \|\u\|_{\bL^{12}}\cdot\|\u\|_{\bW^{1,12/5}}\\
&\leq& C_\Om\|\u\|_{\bW^{1,12/5}}^2\stackrel{(\ref{Int2})}{\leq}
C_\Om(\|A^\frac{5}{8}\u\|^2_{\bL^2}+\|\u\|^2_{\bL^2}),
\de
where the third inequality is due to $\bW^{1,12/5}(\Omega)\subset \bL^{12}(\Om)$.
\end{proof}

\bl\label{Le5}
For any $\frac{3}{4}<\gamma<\beta\leq 1$, there are three positive
continuous functions $F_1, F_3:\mR^2_+\to\mR_+$ and $F_2:\mR_+\to\mR_+$
such that for all $\u,\v\in\bH^{2\beta}$
\ce
&&\|B(\u,\u)-B(\v,\v)\|_{\bL^2}+
\|g^{\nu,\kappa}_N(\|\u\|^2_\infty)\u-g^{\nu,\kappa}_N(\|\v\|^2_\infty)\v\|_{\bL^2}\\
&&\qquad\leq F_1(\|\u\|_{\bH^{2\beta}},\|\v\|_{\bH^{2\beta}})\cdot
\|A^\beta(\u-\v)\|^{\frac{\gamma}{\beta}}_{\bL^2}\cdot
\|\u-\v\|^{\frac{\beta-\gamma}{\beta}}_{\bL^2}\\
&&\qquad\quad+F_2(\|\u\|_{\bH^{2\beta}})\cdot\|A^\beta(\u-\v)\|^{\frac{1}{2\beta}}_{\bL^2}
\cdot\|\u-\v\|^{1-\frac{1}{2\beta}}_{\bL^2}\\
&&\qquad\quad+F_3(\|\u\|_{\bH^{2\beta}},\|\v\|_{\bH^{2\beta}})\cdot
\|\u-\v\|_{\bL^2}.
\de
\el
\begin{proof}
Note that by (iii) of Lemma \ref{Le8}
$$
\|\nabla(\u-\v)\|_{\bL^2}=\|A^{\frac{1}{2}}(\u-\v)\|_{\bL^2}
\leq \|A^\beta(\u-\v)\|^{\frac{1}{2\beta}}_{\bL^2}\cdot\|\u-\v\|^{1-\frac{1}{2\beta}}_{\bL^2}
$$
and
\ce
\|\u-\v\|_\infty&\stackrel{(\ref{PI3})}{\leq}& C_\Om (\|A^\gamma(\u-\v)\|_{\bL^2}
+\|\u-\v\|_{\bL^2})\\
&\leq&  C_\Om(\|A^\beta(\u-\v)\|^{\frac{\gamma}{\beta}}_{\bL^2}\cdot
\|\u-\v\|^{\frac{\beta-\gamma}{\beta}}_{\bL^2}+\|\u-\v\|_{\bL^2}).
\de
The result now follows from
$$
\|B(\u,\u)-B(\v,\v)\|_{\bL^2}\leq \|\u\|_\infty\cdot\|\nabla(\u-\v)\|_{\bL^2}
+\|\u-\v\|_\infty\cdot\|\nabla\v\|_{\bL^2}
$$
and
\ce
\|g^{\nu,\kappa}_N(\|\u\|^2_\infty)\u-g^{\nu,\kappa}_N(\|\v\|^2_\infty)\v\|_{\bL^2}
&\leq&\frac{\kappa}{\nu}\|\u-\v\|_\infty\cdot(\|\u\|_\infty+\|\v\|_\infty)\cdot\|\v\|_{\bL^2}\\
&&+\frac{\kappa}{\nu}\|\u\|^2_\infty\cdot\|\u-\v\|_{\bL^2}.
\de
\end{proof}

We introduce some notations.
Let $I$ be a closed interval of $t$, and let $\mX$ be a Banach space.
By $C(I;\mX)$ we denote the set of all continuous $\mX$-valued functions defined on $I$.
For $0<\theta<1$, $C^\theta(I;\mX)$ means the set of all functions which are strongly
H\"older continuous with the exponent $\theta$. If $I$ is not closed, $v\in C^\theta(I;\mX)$
means that $v\in C^\theta(I_1;\mX)$ for any closed interval $I_1$ contained in $I$.

The following lemma is easily deduced from Lemma \ref{Le8} (cf. \cite{Fu-Ka, Pa}).
\bl\label{Le1}
For $T>0$, let $\f:[0,T]\mapsto\bH^0=\bL^2_\sigma(\Om)$ be continuous and consider
$$
\w(t):=\int^t_0e^{-(t-s)A}\f(s)\dif s.
$$
\begin{enumerate}[(i)]
\item
For any $0\leq\a<\theta<1$
$$
A^\a\w\in C^{1-\theta}([0,T],\bH^0),\quad
\|A^\a\w(t)\|_{\bL^2}\leq C_\a\cdot t^{1-\a}\cdot \sup_{s\in[0,T]}\|\f(s)\|_{\bL^2}.
$$
\item
If $\f\in C([0,T],\bH^0)\cap C^\a((0,T],\bH^0)$ for some $\a\in(0,1)$, then for any $0<\theta<\a$
$$
A\w\in C^\theta((0,T],\bH^0),\quad\p_t\w\in C((0,T],\bH^{2\theta}).
$$
Moreover, $(0,T]$ can be replaced by $[0,T]$ in the above condition and conclusions.
\end{enumerate}
\el
\begin{proof}
The first conclusion is direct from Lemma \ref{Le8}. For the second,
fixing $\delta\in(0,T)$, we write
$$
\w(t)=e^{-(t-\delta)A}\w(\delta)+\int^t_\delta e^{-(t-s)A}\f(s)\dif s=:\Psi_\delta(t)+\Phi_\delta(t).
$$
It is easy to see that
$$
\Psi_\delta(\cdot)\in C^\infty((\delta,T],\bH^0)
$$
and
\ce
\p_t\Phi_\delta(t)&=&e^{-(t-\delta)A}\f(t)-\int^t_\delta
A e^{-(t-s)A}(\f(s)-\f(t))\dif s\\
&=&-A\Phi_\delta(t)+\f(t),\quad \delta\leq t\leq T.
\de
(ii) now follows from Lemma \ref{Le8}.
\end{proof}

For $\a\in[0,1]$, let $\bW^{k+\a,2}(\Om)$ be the complex interpolation
space between $\bW^{k,2}(\Om)$ and $\bW^{k+1,2}(\Om)$.
The following lemma is easily derived by \cite[p.23, Proposition 2.2]{Te} and
the interpolation theorem (cf. \cite{Tr}).
\bl\label{Stoke}
Let $k\in\mN\cup\{0\}$ and $\Om\subset\mR^3$ be a bounded domain of class $C^{k+2}$. For any
$\f\in \bW^{\a,2}(\Om)$, $0\leq\a\leq k$, there exist unique functions
$\u\in\bW^{2+\a,2}(\Om)$ and $P\in W^{1+\a,2}(\Om)$
(with $\int_\cO P\dif x=0$), which solve the following Stokes problem in the distribution sense:
$$
\nu\Delta\u=\nabla P+\f,\ \ \div(\u)=0,\ \
\u|_{\p\Om}=0.
$$
Moreover, there exists a constant $C_{\a,\nu}>0$ such that
$$
\|\u\|_{\bW^{2+\a,2}(\Om)}+\|P\|_{W^{1+\a,2}(\Om)}\leq C_{\a,\nu}\cdot\|\f\|_{\bW^{\a,2}(\Om)}.
$$
\el

\section{Proof of Theorem \ref{Main}}

In this section, we use the following equivalent norm in $\bH^\beta$ ($\beta\in[0,2]$)
$$
\|\u\|_{\bH^\beta}=\|(I+A)^{\beta/2}\u\|_{\bL^2},\ \ \u\in\bH^\beta.
$$
We first prove:
\bl\label{Le01}
For any $\u,\v,\u',\v'\in\bH^2$ we have
\be
&&\<B(\u',\u)-g_N^{\nu,\kappa}(\|\u'\|^2_\infty)\u-B(\v',\v)
+g_N^{\nu,\kappa}(\|\v'\|^2_\infty)\v,(I+A)(\u-\v)\>_{\bL^2}\label{PI1}\\
&&\qquad\leq\frac{3\nu}{4}\|\u-\v\|^2_{\bH^2}+\frac{\nu}{8}\|\u'-\v'\|^2_{\bH^2}
+\frac{\kappa N}{\nu}\cdot \|\w\|^2_{\bH^1}\no\\
&&\qquad\quad+\frac{C_\Omega}{\nu^3}\|\u'-\v'\|^2_{\bH^1}\cdot((1+4\kappa)\|\v\|^2_{\bH^1}+
\kappa\|\u\|^2_{\bH^1})^2.\no
\ee
\el
\begin{proof}
Set
$$
\w=\u-\v,\ \  \ \w'=\u'-\v',
$$
and write (\ref{PI1}) as the following four terms' sum
\ce
I_1&:=&\<B(\u',\w),(I+A)\w\>_{\bL^2},\\
I_2&:=&\<B(\w',\v),(I+A)\w\>_{\bL^2},\\
I_3&:=&-\<g_N^{\nu,\kappa}(\|\u'\|^2_\infty)\w,(I+A)\w\>_{\bL^2},\\
I_4&:=&-\<[g_N^{\nu,\kappa}(\|\u'\|^2_\infty)
-g_N^{\nu,\kappa}(\|\v'\|^2_\infty)]\v,(I+A)\w\>_{\bL^2}.
\de
By $ab\leq \frac{\nu}{4\epsilon} a^2+\frac{\epsilon}{\nu}b^2$, we have
\ce
I_1&\leq&\|B(\u',\w)\|_{\bL^2}\cdot\|(I+A)\w\|_{\bL^2}\\
&\leq&\frac{1}{2\nu}\|(\u'\cdot\nabla)\w\|_{\bL^2}^2+\frac{\nu}{2}\|\w\|_{\bH^2}^2\\
&\leq&\frac{1}{2\nu}\|\u'\|_\infty^2\|\nabla\w\|^2_{\bL^2}+\frac{\nu}{2}\|\w\|_{\bH^2}^2
\de
and similarly,
$$
I_2\leq\frac{1}{\nu}\|\w'\|_\infty^2\|\nabla\v\|^2_{\bL^2}
+\frac{\nu}{4}\|\w\|_{\bH^2}^2.
$$
For $I_3$, by $g^{\nu,\kappa}_N(r)\geq \frac{\kappa}{\nu}(r-N)$ we have
\ce
I_3&=&-g_N^{\nu,\kappa}(\|\u'\|^2_\infty)\cdot\|(I+A)^{1/2}\w\|^2_{\bL^2}\\
&\leq&-\frac{\kappa}{\nu}\|\u'\|^2_\infty\cdot\|\w\|^2_{\bH^1}
+\frac{\kappa N}{\nu}\cdot \|\w\|^2_{\bH^1}.
\de
For $I_4$, by $|g^{\nu,\kappa}_N(r)-g^{\nu,\kappa}_N(r)|\leq\frac{\kappa}{\nu}|r-r'|$ we have
\ce
I_4&\leq&\frac{\kappa}{\nu}\|\w'\|_\infty(\|\u'\|_\infty+\|\v'\|_\infty)
\cdot\|\v\|_{\bH^1}\cdot\|\w\|_{\bH^1}\no\\
&\leq&\frac{\kappa}{\nu}\|\w'\|_\infty\cdot(2\|\u'\|_\infty+\|\w'\|_\infty)
\cdot\|\v\|_{\bH^1}\cdot\|\w\|_{\bH^1}\no\\
&=&\frac{2\kappa}{\nu}(\|\u'\|_\infty\|\w\|_{\bH^1})\cdot(\|\w'\|_\infty\|\v\|_{\bH^1})\no\\
&&+\frac{\kappa}{\nu}\|\w'\|_\infty^2\cdot\|\v\|_{\bH^1}\cdot\|\w\|_{\bH^1}\no\\
&\leq&\frac{\kappa}{2\nu}\|\u'\|^2_\infty\|\w\|_{\bH^1}^2
+\frac{2\kappa}{\nu}\|\w'\|^2_\infty\|\v\|_{\bH^1}^2\no\\
&&+\frac{\kappa}{\nu}\|\w'\|_\infty^2\cdot\|\v\|_{\bH^1}\cdot
(\|\u\|_{\bH^1}+\|\v\|_{\bH^1}).
\de
Combining the above calculations, we obtain
\ce
I_1+I_2+I_3+I_4&\leq&\frac{3\nu}{4}\|\w\|^2_{\bH^2}+\frac{\kappa N}{\nu}\cdot \|\w\|^2_{\bH^1}+
\frac{(1-\kappa)}{2\nu}\|\u'\|^2_\infty\cdot\|\w\|^2_{\bH^1}\\
&&+\frac{1}{\nu}\|\w'\|_\infty^2\cdot((1+3\kappa)\|\v\|^2_{\bH^1}+
\kappa\|\v\|_{\bH^1}\cdot\|\u\|_{\bH^1})\\
&\stackrel{(\ref{PI2})}{\leq}&\frac{3\nu}{4}\|\w\|^2_{\bH^2}+\frac{\kappa N}{\nu}\cdot \|\w\|^2_{\bH^1}+
\frac{(1-\kappa)}{2\nu}\|\u'\|^2_\infty\cdot\|\w\|^2_{\bH^1}\\
&&+\frac{C_\Om}{\nu}\|\w'\|_{\bH^2}\|\nabla\w'\|_{\bL^2}\cdot((1+4\kappa)\|\v\|^2_{\bH^1}+
\kappa\|\u\|^2_{\bH^1})\\
&\stackrel{(\kappa\geq 1)}{\leq}&\frac{3\nu}{4}\|\w\|^2_{\bH^2}+\frac{\nu}{8}\|\w'\|^2_{\bH^2}
+\frac{\kappa N}{\nu}\cdot \|\w\|^2_{\bH^1}\\
&&+\frac{C^2_\Om}{\nu^3}\|\w'\|^2_{\bH^1}\cdot((1+4\kappa)\|\v\|^2_{\bH^1}+
\kappa\|\u\|^2_{\bH^1})^2,
\de
which produces the desired estimate.
\end{proof}

\subsection{Proof of Existence}
Let $\v\in C([0,\infty);\bH^1)\cap L^2_{loc}(\mR_+;\bH^2)$.
We first consider the following linearized equation:
$$
\p_t\u=-\nu A\u+B(\v,\u)-g^{\nu,\kappa}_N(\|\v\|^2_\infty)\u,\ \ \u(0)=\u_0\in\bH^1.
$$
By the standard theory of PDE,
  there is a unique strong solution $\u$ to above equation with
$$
\u\in C([0,\infty);\bH^1)\cap L^2_{loc}(\mR_+;\bH^2).
$$

Let us construct the approximation sequence of equation (\ref{ren}) as follows:
Set $\u_1(t)\equiv0$. For $k=2,3,\cdots$, let
\be
\u_k(t)\in C([0,\infty);\bH^1)\cap L^2_{loc}(\mR_+;\bH^2)\label{PI4}
\ee
solve the following equation
\be
\p_t\u_{k}=-\nu A\u_{k}+B(\u_{k-1},\u_{k})-g^{\nu,\kappa}_N(\|\u_{k-1}\|^2_\infty)\u_{k},\ \ \u_k(0)=\u_0.\label{App}
\ee

Firstly, note that
$$
\<A\u_k,\u_k\>_{\bL^2}=\|\nabla\u_k\|^2_{\bL^2}
$$
and
$$
\<B(\u_{k-1},\u_k),\u_k\>_{\bL^2}=-\<(\u_{k-1}\cdot\nabla)\u_k,\u_k\>_{\bL^2}
=-\frac{1}{2}\<\div\u_{k-1},|\u_k|^2\>_{\bL^2}=0.
$$
By the chain rule, we have from (\ref{App}) that
\be
\dif \|\u_k\|^2_{\bL^2}/\dif t=-2\nu\|\nabla\u_k\|^2_{\bL^2}
-2g^{\nu,\kappa}_N(\|\u_{k-1}\|^2_{\infty})\|\u_k\|^2_{\bL^2}
\leq-2\nu\|\nabla\u_k\|^2_{\bL^2}.\label{Op3}
\ee
Integrating both sides of (\ref{Op3}) yields that
\be
\|\u_k(t)\|^2_{\bL^2}+2\nu\int^t_0\|\nabla\u_k\|^2_{\bL^2}
\dif t\leq \|\u_0\|^2_{\bL^2},\ \ \forall t\geq 0.\label{Es022}
\ee

Secondly, for any $T>0$ we have
\ce
&&\int^T_0|g^{\nu,\kappa}_N(\|\u_{k-1}\|^2_{\infty})|^2\cdot\|\u_k\|^2_{\bL^2}\dif t
\leq \frac{\kappa}{\nu} \sup_{t\in[0,T]}\|\u_k(t)\|^2_{\bL^2}
\int^T_0\|\u_{k-1}\|^4_{\infty}\dif t\\
&&\qquad\stackrel{(\ref{PI2})}{\leq}
C\cdot \sup_{t\in[0,T]}\|\u_k(t)\|^2_{\bL^2}\cdot \sup_{t\in[0,T]}
\|\u_{k-1}(t)\|^2_{\bH^1}\int^T_0\|\u_{k-1}\|^2_{\bH^2}\dif t
\stackrel{(\ref{PI4})}{<}+\infty
\de
and
\be
\int^T_0\|B(\u_{k-1},\u_k)\|^2_{\bL^2}\dif t&\leq&C\int^T_0\|\u_{k-1}\|^2_\infty
\|\nabla\u_k\|^2_{\bL^2}\dif t\no\\
&\stackrel{(\ref{PI3})}{\leq}&C\int^T_0\|\u_{k-1}\|^2_{\bH^2}
\cdot\|\nabla\u_k\|^2\dif t
\stackrel{(\ref{PI4})}{<}+\infty.\label{PP3}
\ee
Thus, recalling $\bH^0=\bL^2_\sigma(\Omega)$, from (\ref{App}) one has
$$
\p_t\u_k\in L^2_{loc}(\mR_+;\bH^0).
$$
Consider the evolution triple
$$
\bH^2\subset\bH^1\subset\bH^0.
$$
By the chain rule (cf. \cite[p.176, Lemma 1.2]{Te}) and Young's inequality, we have
\be
\frac{1}{2}\frac{\dif \|\nabla\u_k\|^2_{\bL^2}}
{\dif t}&=&-\nu\|A\u_k\|^2_{\bL^2}+\<B(\u_{k-1},\u_k),A\u_k\>_{\bL^2}
-g^{\nu,\kappa}_N(\|\u_{k-1}\|^2_{\infty})\|\nabla\u_k\|^2_{\bL^2}\no\\
&\leq&-\frac{\nu}{2}\|A\u_k\|^2_{\bL^2}
+\frac{1}{2\nu}\|\u_{k-1}\|^2_{\infty}\|\nabla\u_k\|^2_{\bL^2}
-g^{\nu,\kappa}_N(\|\u_{k-1}\|^2_{\infty})\|\nabla\u_k\|^2_{\bL^2}\no\\
&\leq&-\frac{\nu}{2}\|A\u_k\|^2_{\bL^2}
+\frac{1-2\kappa}{2\nu}\|\u_{k-1}\|^2_{\infty}\|\nabla\u_k\|^2_{\bL^2}
+\frac{\kappa N}{\nu}\|\nabla\u_k\|^2_{\bL^2},\label{Op11}
\ee
where the last step is due to
$$
g^{\nu,\kappa}_N(r)\geq \frac{\kappa}{\nu}(r-N).
$$
Integrating both sides of (\ref{Op11}) and using (\ref{Es022}) and $\kappa\geq 1$,
we obtain
\be
\|\nabla\u_k(t)\|^2_{\bL^2}+\nu\int^t_0\|A\u_k\|^2_{\bL^2}\dif s
\leq \frac{\kappa N}{\nu^2}\|\u_0\|^2_{\bL^2}+\|\nabla\u_0\|^2_{\bL^2},\ \ \forall t\geq 0.\label{Es02}
\ee

Now set
$$
\w_{k,m}(t):=\u_k(t)-\u_m(t).
$$
Then
\ce
\p_t\w_{k,m}&=&-A\w_{k,m}+B(\u_{k-1},\w_{k,m})+B(\w_{k-1,m-1},\v_m)\\
&&-g^{\nu,\kappa}_N(\|\u_{k-1}\|^2_\infty)\w_{k,m}
-\big[g^{\nu,\kappa}_N(\|\u_{k-1}\|^2_\infty)-g^{\nu,\kappa}_N(\|\v_{m-1}\|^2_\infty)\big]\v_m.
\de
Again, by \cite[p.176, Lemma 1.2]{Te} and Lemma \ref{Le01} we have
\be
\frac{1}{2}\frac{\dif\|\w_{k,m}\|^2_{\bH^1}}{\dif t}
&=&-\nu\|\w_{k,m}\|^2_{\bH^2}+\nu\|\w_{k,m}\|^2_{\bH^1}\no\\
&&+\<B(\u_{k-1},\w_{k,m}),(I+A)\w_{k,m}\>_{\bL^2}\no\\
&&+\<B(\w_{k-1,m-1},\v_m),(I+A)\w_{k,m}\>_{\bL^2}\\
&&-g^{\nu,\kappa}_N(\|\u_{k-1}\|^2_\infty)\<\w_{k,m},(I+A)\w_{k,m}\>_{\bL^2}\no\\
&&-\big[g^{\nu,\kappa}_N(\|\u_{k-1}\|^2_\infty)
-g^{\nu,\kappa}_N(\|\v_{m-1}\|^2_\infty)\big]\<\v_m,\w_{k,m}\>_{\bH^1}\no\\
&\leq&-\frac{\nu}{4}\|\w_{k,m}\|^2_{\bH^2}+
\frac{\nu}{8}\|\w_{k-1,m-1}\|^2_{\bH^2}+(\nu+\frac{\kappa N}{\nu})\|\w_{k,m}\|^2_{\bH^1}\no\\
&&+\frac{C_\Omega}{\nu^3}\|\w_{k-1,m-1}\|^2_{\bH^1}\cdot((1+4\kappa)\|\u_{k}\|^2_{\bH^1}+
\kappa\|\u_m\|^2_{\bH^1})^2.\label{Lp2}
\ee
Integrating this inequality and using (\ref{Es022}) and (\ref{Es02}), we get
\ce
&&\|\w_{k,m}(t)\|^2_{\bH^1}+\frac{\nu}{2}\int^t_0\|\w_{k,m}\|^2_{\bH^2}\dif s\\
&&\qquad \leq\frac{\nu}{4}\int^t_0\|\w_{k-1,m-1}\|^2_{\bH^2}\dif s
+C_{\nu,N}\int^t_0\|\w_{k,m}\|^2_{\bH^1}\dif s\\
&&\qquad \quad+C_{\nu,N,\|\u_0\|_{\bH^1}}\cdot\int^t_0\|\w_{k-1,m-1}\|^2_{\bH^1}\dif s.
\de

Set
$$
h(t):=\varlimsup_{k,m\to\infty}\sup_{s\in[0,t]}\|\w_{k,m}(s)\|^2_{\bH^1}
$$
and
$$
f(t):=\varlimsup_{k,m\to\infty}\int^t_0\|\w_{k,m}\|^2_{\bH^2}\dif s.
$$
Then by (\ref{Es022}), (\ref{Es02}) and Fatou's lemma, we have
$$
\frac{\nu}{2}f(t)\leq\frac{\nu}{4}f(t)+C_{\nu,N,\|\u_0\|_{\bH^1}}\int^t_0h(s)\dif s
$$
and
$$
h(t)\leq \frac{\nu}{4}f(t)+C_{\nu,N,\|\u_0\|_{\bH^1}}\int^t_0h(s)\dif s\leq
2C_{\nu,N,\|\u_0\|_{\bH^1}}\int^t_0g(s)\dif s.
$$
Hence, by Gronwall's inequality we have
$$
h(t)=f(t)=0,\ \ \forall t\geq 0.
$$
Thus, there exists a function $\u\in C([0,\infty);\bH^1)\cap L^2_{loc}(\mR_+;\bH^2)$ such that
for any $T>0$
\ce
\varlimsup_{k\to\infty}\sup_{s\in[0,T]}\|\u_{k}-\u\|^2_{\bH^1}
+\varlimsup_{k\to\infty}\int^T_0\|\u_{k}(s)-\u(s)\|^2_{\bH^2}\dif s=0.
\de
Lastly, taking limits $k\to\infty$ for
$$
\u_{k}(t)=\u_0-\int^t_0A\u_{k}\dif s+\int^t_0B(\u_{k-1},\u_{k})\dif s-
\int^t_0g(\|\u_{k-1}\|^2_\infty)\u_{k}\dif s
$$
and inequalities (\ref{Es022}) and (\ref{Es02}),
we can see that $\u(t)$ satisfies (\ref{Op1}),  (\ref{Es22}) and (\ref{Es222}).

\subsection{Proof of Decay Estimate (\ref{decay})}

Following the method of Heywood \cite{He}, by the chain rule
and \cite[p.649 (14)]{He}, we have
\ce
\frac{\dif \|\nabla\u\|^2_{\bL^2}}
{\dif t}&=&-2\nu\|A\u\|^2_{\bL^2}+2\<B(\u,\u),A\u\>_{\bL^2}
-2g^{\nu,\kappa}_N(\|\u\|^2_{\infty})\|\nabla\u\|^2_{\bL^2}\no\\
&=&-\nu\|A\u\|^2_{\bL^2}+C_{\nu,\Omega}\|\nabla\u\|^4_{\bL^2}
+C'_{\nu,\Om}\|\nabla\u\|^6_{\bL^2}.
\de
Note that
$$
\int^\infty_0\|\nabla\u_k\|^2_{\bL^2}
\dif t\leq \frac{\|\u_0\|^2_{\bL^2}}{2\nu}=:\Lambda.
$$
In Lemma \ref{Le2}, if we take $\phi(t)=\|\nabla\u(t)\|^2_{\bL^2}$,
$\beta=1/(C'_{\nu,\Om}\Lambda)$ and $\a=C_{\nu,\Omega}+1/\Lambda$,
then
\ce
\|\nabla\u(t)\|^2_{\bL^2}\leq \frac{e^{C_{\nu,\Omega}\Lambda+1}-1}
{(C_{\nu,\Omega}+1/\Lambda)t}
\leq C_3 t^{-1}
\de
for $t\geq (C'_{\nu,\Om}\Lambda)e^{C_{\nu,\Omega}\Lambda+1}=T^*$. Thus, (\ref{decay})
follows.

\subsection{Proof of Continuous Dependence (\ref{Col})}

Set
$$
\w_{N,M}(t):=\u_N(t)-\v_M(t).
$$
Once again, by the chain rule (cf. \cite[p.176, Lemma 1.2]{Te})  we have
\ce
\frac{1}{2}\frac{\dif\|\w_{N,M}\|^2_{\bH^1}}{\dif t}
&=&-\|\w_{N,M}\|^2_{\bH^2}+\|\w_{N,M}\|^2_{\bH^1}
+\<B(\u_N,\w_{N,M}),\w_{N,M}\>_{\bH^1}\\
&&+\<B(\w_{N,M},\v_M),\w_{N,M}\>_{\bH^1}
-g^{\nu,\kappa}_N(\|\u_N\|^2_\infty)\|\w_{N,M}\|^2_{\bH^1}\\
&&-\big[g^{\nu,\kappa}_N(\|\u_N\|^2_\infty)
-g^{\nu,\kappa}_N(\|\v_M\|^2_\infty)\big]\<\v_M, \w_{N,M}\>_{\bH^1}\\
&&-\big[g^{\nu,\kappa}_N(\|\v_M\|^2_\infty)
-g^{\nu,\kappa}_M(\|\v_M\|^2_\infty)\big]\<\v_M, \w_{N,M}\>_{\bH^1}.
\de
Noting that
$$
|g^{\nu,\kappa}_N(r)-g_M(r)|\leq \frac{\kappa}{\nu}|N-M|,\ \ \forall r\geq 0,\ \ N,M\geq 1,
$$
as in the proof of existence, by Lemma \ref{Le01} and Young's inequality we find that
\ce
\frac{1}{2}\frac{\dif\|\w_{N,M}\|^2_{\bH^1}}{\dif t}
&\leq&-\frac{\nu}{8}\|\w_{N,M}\|^2_{\bH^2}+(\nu+\frac{\kappa N}{\nu})\|\w_{N,M}\|^2_{\bH^1}\\
&&+\frac{C_\Omega}{\nu^3}\|\w_{N,M}\|^2_{\bH^1}\cdot((1+4\kappa)\|\u_{N}\|^2_{\bH^1}+
\kappa\|\v_M\|^2_{\bH^1})^2\\
&&+\frac{\kappa}{\nu}|N-M|\cdot\|\v_M\|_{\bL^2}\cdot\|\w_{N,M}\|_{\bH^2}\\
&\stackrel{(\ref{Es22})(\ref{Es222})}{\leq}&
 -\frac{\nu}{16}\|\w_{N,M}\|^2_{\bH^2}
+C_{\nu,N,M,\u_0,\v_0}\|\w_{N,M}\|^2_{\bH^1}+C_{\nu,\v_0} |N-M|^2,
\de
where $C_{\nu,\v_0}=4\kappa^2\|\v_0\|^2_{\bL^2}/\nu^3$,
$$
C_{\nu,N,M,\u_0,\v_0}=(\nu+\frac{\kappa N}{\nu})
+\frac{C_\Omega}{\nu^3}((1+4\kappa)K_{\nu,N,\u_0}+
\kappa K_{\nu,M,\v_0})^2
$$
and
$$
K_{\nu,N,\u_0}:=\frac{\kappa N}{\nu^2}\|\u_0\|^2_{\bL^2}+\|\u_0\|^2_{\bH^1}.
$$
The estimate (\ref{Col}) now follows by Gronwall's inequality.

\section{Proof of Theorem \ref{main}}

\subsection{Proof of Part (I)}
Let $\u(t)$ be the unique strong solution of equation (\ref{ren}).
By Duhamel's formula, we may write
\be
\u(t)&=&e^{-A t}\u_0+\int^t_0 e^{-(t-s)A}B(\u,\u)\dif s-\int^t_0e^{-(t-s)A}
g^{\nu,\kappa}_N(\|\u\|^2_\infty)\u\dif s\no\\
&=:&\w_1(t)+\w_2(t)+\w_3(t).\label{Lp1}
\ee

First of all, it is clear that $\w_1\in C^\infty((0,T];\bH^2)$ and
\be
\|A\w_1(t)\|_{\bL^2}\leq \|A\u_0\|_{\bL^2}.\label{Op8}
\ee
For $\w_2(t)$, by (i) of Lemma \ref{Le8} we have
\be
\|A^{\frac{5}{8}} \w_2(t)\|_{\bL^2}&\leq&\int^t_0
\|A^{\frac{7}{8}}e^{-(t-s)A}A^{-\frac{1}{4}}B(\u,\u)\|_{\bL^2}\dif s\no\\
&\leq&\int^t_0\frac{\|A^{-\frac{1}{4}}B(\u,\u)\|^2_{\bL^2}}{(t-s)^{\frac{7}{8}}}\dif s
\stackrel{(\ref{PI2})}{\leq} C_\Om\cdot\int^t_0\frac{\|A^\frac{1}{2}
\u\|^2_{\bL^2}}{(t-s)^{\frac{7}{8}}}\dif s\no\\
&\stackrel{(\ref{Es222})}{\leq}& C_\Om   K_{\nu,N,\u_0}\cdot t^{\frac{1}{8}},\label{Op9}
\ee
where
$$
K_{\nu,N,\u_0}:=\frac{\kappa N}{\nu^2}\|\u_0\|^2 +\|\nabla\u_0\|^2.
$$
For $\w_3(t)$, recalling (\ref{Id}) and
by Lemma \ref{Le8},  we have for $\a\in[1/2,1)$
\ce
\|A^\a\w_3(t)\|_{\bL^2}&\leq&\int^t_0g^{\nu,\kappa}_N(\|\u\|^2_\infty)\cdot
\|A^\a e^{-(t-s)A}\u\|_{\bL^2}\dif s\no\\
&\leq&\int^t_0\|\u\|^2_\infty\cdot \frac{1}{(t-s)^{\a-\frac{1}{2}}}\cdot
\|A^{\frac{1}{2}}\u\|_{\bL^2}\dif s\no\\
&\stackrel{(\ref{PI2})}{\leq}&C_\Om\int^t_0\|\u\|_{\bH^2}\cdot\|\nabla\u\|^2_{\bL^2}\cdot
\frac{1}{(t-s)^{\a-\frac{1}{2}}}\dif s\no\\
&\leq&C_\Om\int^t_0\|A\u\|_{\bL^2}\cdot\|\nabla\u\|^2_{\bL^2}\cdot
\frac{1}{(t-s)^{\a-\frac{1}{2}}}\dif s\no\\
&&+C_\Om\int^t_0\|\u\|_{\bL^2}\cdot\|\nabla\u\|^2_{\bL^2}\cdot
\frac{1}{(t-s)^{\a-\frac{1}{2}}}\dif s\no\\
&=&:I_1+I_2.
\de
By (\ref{Es22}), (\ref{Es222}) and H\"older's inequality we have
$$
I_1\leq C_\Om  K_{\nu,N,\u_0} t^{1-\a}
\left(\int^t_0\|A\u(s)\|^2_{\bL^2}\dif s\right)^{1/2}
\leq C_\Om  K_{\nu,N,\u_0}^{3/2}\cdot  t^{1-\a}
$$
and
$$
I_2\leq C_\Om\|\u_0\|_{\bL^2} K_{\nu,N,\u_0}^{\frac{1}{2}}  t^{1-\a}
 \left(\int^t_0\|\nabla\u(s)\|^2_{\bL^2}\dif s\right)^{1/2}
\leq C_\Om\|\u_0\|^2_{\bL^2} K_{\nu,N,\u_0}^{1/2}\cdot  t^{1-\a}.
$$
Hence
\be
\|A^\a\w_3(t)\|_{\bL^2}\leq C_\Om\Big(\|\u_0\|^2_{\bL^2} K_{\nu,N,\u_0}^{1/2}
+ K_{\nu,N,\u_0}^{3/2}\Big)\cdot t^{1-\a}.\label{Op10}
\ee

Combining  (\ref{Lp1}), (\ref{Op8}), (\ref{Op9}) and (\ref{Op10}), we find that
\ce
\|A^{\frac{5}{8}}\u(t)\|_{\bL^2}&\leq& \|A^{\frac{5}{8}}\u_0\|_{\bL^2}
+C_\Om   K_{\nu,N,\u_0}\cdot t^{\frac{1}{8}}\\
&&+C_\Om\Big(\|\u_0\|^2_{\bL^2} K_{\nu,N,\u_0}^{1/2}
+ K_{\nu,N,\u_0}^{3/2}\Big)\cdot t^{3/8}\\
&=:&M_0(t,\nu,N,\u_0)
\de
and by (\ref{PI6}) and (\ref{Es22})
\ce
\|B(\u(t),\u(t))\|_{\bL^2}&\leq& C_\Om(\|A^{\frac{5}{8}}\u(t)\|^2_{\bL^2}+\|\u(t)\|^2_{\bL^2})\\
&\leq& C_\Om\cdot (M_0(t,N,\u_0)^2+\|\u_0\|^2_{\bL^2}).
\de
By  (1$^0$) of Lemma \ref{Le1}, we have for any $\frac{3}{4}<\gamma<1$
\be
\|A^\gamma\u(t)\|_{\bL^2}&\leq& \|A^\gamma\u_0\|_{\bL^2}+
C_\Om\cdot (M_0(t,\nu,N,\u_0)^2+\|\u_0\|^2_{\bL^2})\cdot t^{1-\gamma}\no\\
&&+C_\Om\Big(\|\u_0\|^2_{\bL^2} K_{\nu,N,\u_0}^{1/2}
+ K_{\nu,N,\u_0}^{3/2}\Big)\cdot t^{1-\gamma},\label{Op12}
\ee
which then yields the estimate (\ref{Es1}) by  (\ref{PI3}).

\subsection{Proof of Part (II)}
In this subsection, we assume $\Om=\mR^3$ or $\Om$ is a bounded uniform $C^4$-domain.
Our proof is concentrated on the case of bounded domain. Clearly, it also works for $\Om=\mR^3$.

Below, fix $T>0$ and set
$$
\f(s):=B(\u(s),\u(s))-g^{\nu,\kappa}_N(\|\u(s)\|^2_\infty)\u(s).
$$
Then by Lemma \ref{Le5} and (\ref{PI3}), (\ref{Op12})
$$
[0,T]\ni s\mapsto\f(s)\in\bH^0=\bL^2_\sigma(\Om)\ \mbox{ is continuous.}
$$
By (i) of Lemma \ref{Le1}, we have for any $\beta\in(0,1)$ and $0<\theta<1$
\be
\u\in C^\theta([0,T],\bH^0)\cap C([0,T];\bH^{2\beta}).\label{Po5}
\ee
Thus, by Lemma \ref{Le5} and (\ref{PI3}), for any $\frac{3}{4}<\gamma<\beta\leq 1$,
there are constants $C_1, C_2,C_3>0$
such that for all $t,s\in[0,T]$
\be
\|\f(t)-\f(s)\|_{\bL^2}&\leq& C_1\cdot
\|A^\beta(\u(t)-\u(s))\|^{\frac{\gamma}{\beta}}_{\bL^2}\cdot
\|\u(t)-\u(s)\|^{1-\frac{\gamma}{\beta}}_{\bL^2}\no\\
&&+C_2\cdot\|A^\beta(\u(t)-\u(s))\|^{\frac{1}{2\beta}}_{\bL^2}
\cdot\|\u(t)-\u(s)\|^{1-\frac{1}{2\beta}}_{\bL^2}\no\\
&&+C_3\cdot\|\u(t)-\u(s)\|_{\bL^2}.\label{PP4}
\ee
Choosing $\beta$ close to $1$ and $\gamma$ close to $\frac{3}{4}$
and using (\ref{Op12}) and (\ref{Po5}), we  find that for any $0<\a<\frac{1}{4}$
$$
\f\in C^\a([0,T],\bH^0).
$$
Thus, by (ii) of Lemma \ref{Le1} and (\ref{Lp1}) we have, for any $0<\a<\frac{1}{4}$
$$
A\u\in C^\a((0,T],\bH^0),\quad \p_t\u\in C((0,T],\bH^{2\a})
$$
Using induction and (\ref{PP4}) with $\beta=1$ as well as (\ref{Po5}), one finds that
for any $n\in\mN$ and $0<\a<1-\left(\frac{3}{4}\right)^{n+1}$
$$
\f\in C^\a((0,T],\bH^0)
$$
and
$$
A\u\in C^\a((0,T],\bH^0),\ \ \p_t\u\in C((0,T],\bH^{2\a}).
$$
In particular, by $\bH^\a\subset\bW^{\a,2}(\Om)$ for $\a\in[0,2]$ we have
\be\label{Pl1}
\u\in C^\a((0,T],\bW^{2,2}(\Om)),\ \ \p_t\u\in C((0,T],\bW^{9/5,2}(\Om)).
\ee

Set
$$
\b(t):=(\u(t)\cdot\nabla)\u(t)+g^{\nu,\kappa}_N(\|\u(t)\|^2_\infty)\u(t).
$$
As in the proof of Lemma \ref{Le5}, it is not hard to verify by (\ref{Pl1}) that
\be
\b(t)\in C((0,T];\bW^{1,2}(\Om)).\label{ES22}
\ee
Consider the Stokes equation:
\ce
\left\{
\begin{aligned}
&\nu\Delta\u+\nabla P=\p_t\u+\b\ \mbox{ in $\Om$},\\
&\dif\u=0\ \mbox{ in $\Om$}, \ \u|_{\p\Om}=0.
\end{aligned}
\right.
\de
By (\ref{Pl1}), (\ref{ES22}) and Lemma \ref{Stoke} with $\a=1$, we have
$$
\u(t)\in C((0,T],\bW^{3,2}(\Om)).
$$
As above, a simple calculation shows that
\be
\b(t)\in C((0,T];\bW^{2,2}(\Om)). \label{ES33}
\ee
By (\ref{Pl1}), (\ref{ES33}) and Lemma \ref{Stoke} with $\a=\frac{9}{5}$ again, we further have
\be
\u(t)\in  C((0,T],\bW^{\frac{19}{5},2}(\Om))\label{Eo00}
\ee
and
\be
P\in C((0,T],W^{\frac{14}{5},2}(\Om)).\label{Eo}
\ee
By (\ref{Pl1}), (\ref{Eo00}), (\ref{Eo}) and
the Sobolev embedding theorem (cf. \cite[Theorem 4.6.1]{Tr}),
we finally obtain that
$$
\p_t\u,~\p_i\u,~\p_i\p_j\u,~ P
\in C((0,T]\times\bar\Om;\mR^3)
$$
and (\ref{Cl}) holds.

\subsection{Proof of Part (III)}

In this subsection, we assume $\Om=\mR^3$.
\bl\label{Lem1}
For fixed $q\geq 2$ and $r\geq 1$, there exists $\kappa:=\kappa(q):=Cq^4$,
where $C$ is a universal constant, such that for any $N\geq 1$ and $t\geq0$
\be
\|\u^\kappa_N(t)\|^{r}_{\bL^q}\leq \|\u_0\|^{r}_{\bL^q}
+\frac{r\kappa}{\nu}\cdot N\int^t_0\|\u^\kappa_N\|^{r}_{\bL^q}\dif s\label{PP1}
\ee
and
\be
\int^t_0\|\u^\kappa_N\|^2_\infty\cdot\|\u^\kappa_N\|^{r}_{\bL^q}\dif s
\leq\frac{2\nu}{r\kappa}\|\u_0\|^{r}_{\bL^q}
+2N\int^t_0\|\u^\kappa_N\|^{r}_{\bL^q}\dif s.\label{PP0}
\ee
\el
\begin{proof}
Let $\u:=\u^\kappa_N$.
Taking the scalar product for both sides of equation (\ref{Cl}) with $q|\u|^{q-2}\u$,
and then integrating over $\mR^3$, we find by the integration by parts formula
\ce
\frac{\dif\|\u\|^{q}_{\bL^q}}{\dif t}&=&-q\nu\||\nabla\u||\u|^{(q-2)/2}\|^2_{\bL^2}
-\frac{4(q-2)\nu}{q}\|\nabla|\u|^{q/2}\|^2_{\bL^2}\no\\
&&+q\<\nabla P, |\u|^{q-2}\u\>_{\bL^2}
-q\cdot g^{\nu,\kappa}_N(\|\u\|^2_\infty)\|\u\|^q_{\bL^q},
\de
where we have used that
$$
q\<(\u\cdot\nabla)\u, |\u|^{q-2}\u\>_{\bL^2}=\<\u,\nabla|\u|^{q/2}\>_{\bL^2}=0.
$$
Let $f$ be an increasing smooth function on $[0,\infty)$.
We further have
\be
\frac{\dif f(\|\u\|^{q}_{\bL^q})}{\dif t}&=&f'(\|\u\|^{q}_{\bL^q})
\Big[-q\nu\||\nabla\u||\u|^{(q-2)/2}\|^2_{\bL^2}
-\frac{4(q-2)\nu}{q}\|\nabla|\u|^{q/2}\|^2_{\bL^2}\no\\
&&+q\<\nabla P, |\u|^{q-2}\u\>_{\bL^2}
-q\cdot g^{\nu,\kappa}_N(\|\u\|^2_\infty)\|\u\|^q_{\bL^q}\Big].\label{Es2}
\ee

On the other hand, taking the divergence for equation (\ref{Cl}) we have
$$
\Delta P=\div[(\u\cdot\nabla)\u],
$$
which gives
$$
P=-(-\Delta)^{-1}\div[(\u\cdot\nabla\u)]
=-(-\Delta)^{-1}\p_j\p_i(u^i\cdot u^j).
$$
So, by the Calder\'on-Zygmund inequality we get for any $\gamma\geq 2$ (cf. \cite{St})
\be
\|P\|_{L^\gamma}\leq C_1\cdot\gamma\cdot\|\u\|_{\bL^{2\gamma}}^2,\label{Lp3}
\ee
Here and below, $C_i, i=1,2,3$ are universal constants.
Thus, by Young's inequality and H\"older's inequality we have
\ce
q\<\nabla P, |\u|^{q-2}\u\>_{\bL^2}&=&q\<P, |\u|^{q-2}\div\u\>_{\bL^2}
+q\<P, \nabla|\u|^{q-2}\cdot \u\>_{\bL^2}\\
&\leq&q\nu\||\nabla\u|\cdot |\u|^{(q-2)/2}\|^2_{\bL^2}+\frac{C_2 q^3}{\nu}\||P|\cdot|\u|^{(q-2)/2}\|^2_{L^2}\\
&\leq&q\nu\||\nabla\u|\cdot |\u|^{(q-2)/2}\|^2_{\bL^2}+\frac{C_2 q^3}{\nu}\|\u\|^{q-2}_{\bL^{q+2}}\cdot\|P\|^2_{L^{(q+2)/2}}\\
&\stackrel{(\ref{Lp3})}{\leq}&
q\nu\||\nabla\u|\cdot |\u|^{(q-2)/2}\|^2_{\bL^2}+\frac{C_3 q^5}{\nu}\|\u\|^{q+2}_{\bL^{q+2}}\\
&\leq&q\nu\||\nabla\u|\cdot |\u|^{(q-2)/2}\|^2_{\bL^2}+\frac{C_3 q^5}{\nu}\|\u\|^2_\infty\cdot\|\u\|^{q}_{\bL^q}.
\de
Substituting this estimate into (\ref{Es2}), we get
\ce
&&\frac{\dif f(\|\u\|^{q}_{\bL^q})}{\dif t}+\frac{4(q-2)\nu}{q}f'(\|\u\|^{q}_{\bL^q})\cdot
\|\nabla|\u|^{q/2}\|^2_{\bL^2}\\
&&\qquad \leq f'(\|\u\|^{q}_{\bL^q})\cdot\Big[\frac{C_3 q^5}{\nu}\cdot\|\u\|^2_\infty\cdot\|\u\|^{q}_{\bL^q}
-q\cdot g^{\nu,\kappa}_N(\|\u\|^2_\infty)\|\u\|^q_{\bL^q}\Big].
\de

Now noticing that
$$
g^{\nu,\kappa}_N(r)\geq \frac{\kappa(r-N)}{\nu}, \ \ r\geq 0,
$$
we find that if
\be
\kappa=2C_3 q^4,\label{Ka}
\ee
then
\be
&&\frac{\dif f(\|\u\|^{q}_{\bL^q})}{\dif t}+\frac{4(q-2)\nu}{q}f'(\|\u\|^{q}_{\bL^q})\cdot
\|\nabla|\u|^{q/2}\|^2_{\bL^2}\no\\
&&\qquad \leq f'(\|\u\|^{q}_{\bL^q})\cdot\Big[\frac{q\kappa}{\nu}\cdot N\cdot\|\u\|^{q}_{\bL^q}
-\frac{q\kappa}{2\nu}\cdot \|\u\|^2_\infty\cdot\|\u\|^q_{\bL^q}\Big].
\label{OP1}
\ee

Lastly, taking $f_\epsilon(x):=(\epsilon+x)^{r/q}$ in (\ref{OP1}), then integrating with respect to $t$
and letting $\epsilon\downarrow 0$  yield (\ref{PP1}) and (\ref{PP0}).
\end{proof}

\bl
Fix  $r_0\geq 1$ and $q_0\geq 2$. Let $\u_0\in\bH^2_2$ and set $N_0:=C\|\u_0\|^2_{\bH^2_2}$
for some universal constant $C$.
There exists  $n_0:=n_0(\nu,N_0,q_0,r_0)$ large enough such that
for all $n\geq n_0$, $N\geq N_0\vee 1$ and $t\geq 0$
\be
\|\u_N^{\kappa_n}(t)\|_{\bL^{q_0+2q_0 n/r_0}}\leq 2N^{\frac{n}{2n+r_0}}+2N^{\frac{n+1}{2n+r_0}}
\left[\int^t_0\|\u^{\kappa_n}_N\|^{r_0}_{\bL^{q_0}}\dif s\right]^{\frac{1}{2n+r_0}},\label{Ep1}
\ee
where $\kappa_n=2C_3\cdot [(2n+r_0)q_0/r_0]^4$ (see (\ref{Ka})).
In particular, there is an $n_0:=n_0(\nu,N_0)$ large enough such that
for all $n\geq n_0$, $N\geq N_0\vee 1$ and $t\geq 0$
\be
\|\u_N^{\kappa_n}(t)\|_{\bL^{6(n+1)}}\leq 3N^{\frac{1}{2}}.\label{Ep2}
\ee
\el
\begin{proof} Let $\u:=\u^\kappa_N$.
First of all, by the Gagliado-Nireberg  inequality (\ref{GN}),
there is a universal constant $C_0\geq 1$ such that for any $q\geq 2$
\be
\|\u_0\|_{\bL^q}\leq \|\u_0\|^{1-2/q}_{\infty}\cdot\|\u_0\|^{2/q}_{\bL^2}
\leq C_0\|\u_0\|_{\bH^2_2}=:N^{1/2}_0.\label{OO1}
\ee

Define
$$
q_n:=q_{n-1}+2q_0/r_0=(2n+r_0)q_0/r_0\ \
$$
and
$$
r_n:=r_0 q_n/q_0=2n+r_0.
$$
Then we have
\ce
\int^t_0\|\u\|_{\bL^{q_{n+1}}}^{r_{n+1}}\dif s
&\leq&\int^t_0\|\u\|_{\infty}^2\cdot\|\u\|_{\bL^{q_{n}}}^{r_{n}}\dif s\\
\mbox{(by (\ref{PP0})) }&\leq&\frac{2\nu }{r_n\cdot\kappa(q_n)}\|\u_0\|^{r_n}_{\bL^{q_n}}
+2N\int^t_0\|\u\|_{\bL^{q_n}}^{r_n}\dif s\\
\mbox{(by $r_n\kappa(q_n)\geq 2$) }
&\leq&\nu\|\u_0\|^{r_n}_{\bL^{q_n}}
+2N\int^t_0\|\u\|_{\bL^{q_n}}^{r_n}\dif s\\
\mbox{(by iterating) }
&\leq&\nu\sum_{k=0}^n(2N)^k\|\u_0\|^{r_{n-k}}_{\bL^{q_{n-k}}}
+(2N)^{n+1}\int^t_0\|\u\|_{\bL^{q_0}}^{r_0}\dif s\\
\mbox{(by (\ref{OO1}) )}
&\leq&\nu\sum_{k=0}^n(2N)^k\cdot N_0^{r_{n-k}/2}
+(2N)^{n+1}\int^t_0\|\u\|_{\bL^{q_0}}^{r_0}\dif s\\
\mbox{(by $r_n=2n+r_0$)}
&=&\frac{\nu((2N)^{n+1}N_0^{r_0/2}-N_0^{n+1+r_0/2})}{2N-N_0}\\
&&+(2N)^{n+1}\int^t_0\|\u\|_{\bL^{q_0}}^{r_0}\dif s\\
\mbox{(by $N\geq N_0\vee 1$)}&\leq&\frac{\nu N_0^{r_0/2}(2N)^{n+1}}{N}
+(2N)^{n+1}\int^t_0\|\u\|_{\bL^{q_0}}^{r_0}\dif s.
\de
Hence, by (\ref{PP1})
\ce
\|\u(t)\|^{r_n}_{\bL^{q_n}}&\leq& \|\u_0\|_{\bL^{q_n}}^{r_n}
+r_n\kappa(q_n)\cdot N_0^{r_0/2}(2N)^{n}\\
&&+\frac{r_n\kappa(q_n)N}{\nu}\cdot(2N)^{n}\int^t_0\|\u\|^{r_0}_{\bL^{q_0}}\dif s.
\de
Now taking the root $1/r_n$ and noting that
$$
\lim_{n\rightarrow\infty}(r_n\kappa(q_n))^{1/r_n}
\stackrel{(\ref{Ka})}{=}\lim_{n\rightarrow\infty}(2C_3r_n q^4_n)^{1/r_n}=1,
$$
we obtain the desired estimate (\ref{Ep1}).

As for (\ref{Ep2}), it follows by taking $r_0=2$ and $q_0=6$ in (\ref{Ep1}) and
noting that
$$
\int^t_0\|\u\|^2_{\bL^6}\dif s\stackrel{(\ref{Int})}{\leq}
C\int^t_0\|\nabla\u\|^2_{\bL^2}\dif s\stackrel{(\ref{Es22})}{\leq} \frac{C}{2\nu}\|\u_0\|^2_{\bL^2}.
$$
The proof is complete.
\end{proof}

We are now in a position to give

{\bf Proof of (\ref{Es}):} By (\ref{Int2}), (\ref{Op12}) and  (\ref{Es22}) we have
\be
\|\u^\kappa_N(t)\|_{\bW^{1,4}}&\leq& C(\|A^{7/8}\u^\kappa_N(t)\|_{\bL^2}+\|\u^\kappa_N(t)\|_{\bL^2})\\
&\leq& C\|A^{7/8}\u_0\|_{\bL^2}+
C\cdot (M_0(t,\nu,N,\u_0)^2+\|\u_0\|^2_{\bL^2})\cdot t^{1/8}\no\\
&&+C_\Om\Big(\|\u_0\|^2_{\bL^2} K_{\nu,N,\u_0}^{1/2}
+ K_{\nu,N,\u_0}^{3/2}\Big)\cdot t^{1/8}+C\|\u_0\|_{\bL^2}.\label{Po6}
\ee
By the  Gagliado-Nireberg  inequality (\ref{GN}) and (\ref{Ep2}), we have
$$
\|\u^\kappa_N(t)\|_{\infty}\leq C\|\u^\kappa_N(t)\|_{\bW^{1,4}}^{\frac{2}{n+3}}
\cdot\|\u^\kappa_N(t)\|_{\bL^{6(n+1)}}^{\frac{n+1}{n+3}}.
$$
Letting $n$ be large enough, the estimate (\ref{Es}) follows from (\ref{Ep2}) and (\ref{Po6}).

\section{Proof of Theorem \ref{Th3}}

We need the following simple lemma. For the reader's convenience, a short proof is provided here.
\bl\label{Le14}
Let $(\mX,\|\cdot\|_\mX)$ be a uniformly convex
Banach space and $K\subset\mX$. Then $K$ is relatively compact in $\mX$
if and only if there exists a family of finite dimensional subspaces $\{\mX_n,n\in\mN\}$ of $\mX$
such that
\be
\sup_{n\in\mN}\sup_{x\in K}\|\Pi_n x\|_{\mX}<+\infty\label{LP1}
\ee
and
\be
\lim_{n\rightarrow\infty}\sup_{x\in K}\|(I-\Pi_n)x\|_{\mX}=0,\label{LP2}
\ee
where $\Pi_n$ is the projection operator from $\mX$ to $\mX_n$, i.e., $\Pi_nx\in\mX_n$
is the unique element such that
$$
\|x-\Pi_nx\|_{\mX}=\inf_{y\in\mX_n}\|x-y\|_{\mX}.
$$
\el
\begin{proof}
(``Only if'':) Let $K$ be relatively compact in $\mX$. For any $n\in\mN$, there are finite points
$\{x_1,\cdots,x_m\}\subset K$ such that
$$
K\subset\cup_{i=1}^m B_{1/n}(x_i),
$$
where $B_{1/n}(x_k)$ denotes the ball in $\mX$ with center $x_k$ and radius $1/n$. Now put
$$
\mX_n:=\mbox{span}\{x_1,\cdots,x_m\}.
$$
It is easy to see that the corresponding $\Pi_n$ satisfy (\ref{LP1}) and (\ref{LP2}).

(``If'':) Fix any sequence $\{x_k,k\in\mN\}\subset K$. It suffices to prove that there is a subsequence
$x_{k_l}$ such that  $x_{k_l}$ converges to some point $x\in\mX$.
For any $n\in\mN$, since $\mX_n$ is finite dimensional, by (\ref{LP1})
there is a subsequence $x_{k^{(n)}_l}$ and $y_n\in\mX_n$ such that $\Pi_nx_{k^{(n)}_l}$ converges to $y_n$
as $l\rightarrow\infty$. By the diagonalization method, one can find a common subsequence $x_{k_l}$
such that for any $n\in\mN$
$$
\lim_{l\rightarrow\infty}\|\Pi_nx_{k_l}-y_n\|_\mX=0.
$$
Noting that
$$
\|y_n-y_m\|_{\mX}\leq \|\Pi_nx_{k_l}-y_n\|_\mX+\|\Pi_nx_{k_l}-y_n\|_\mX+\|\Pi_nx_{k_l}-P_my_n\|_\mX,
$$
we have by (\ref{LP2}) that $\{y_n,n\in\mN\}$ is a Cauchy sequence in $\mX$. So, there is an $x\in\mX$
such that $y_n$ converges to $x$ in $\mX$. By (\ref{LP2}) again, it is easy to find that
$x_{k_l}$ converges to $x$ in $\mX$. The proof is complete.
\end{proof}

Since we have assumed that $\Om$ is a bounded domain in Theorem \ref{Th3},
$\bH^1=\bW^{1,2}_{0,\sigma}(\Om)$ is compactly embedded in $\bH^0=\bL^2_\sigma(\Om)$.
Let $0<\l_1\leq\l_2\leq\cdots\leq\l_k\to\infty$ be the eigenvalues of $A$, and
$\sE:=\{\e_k; k\in\mN\}$ the corresponding orthonormal eigenvectors, i.e.,
\be
A\e_k=\l_k\e_k,\ \ \<\e_k,\e_j\>_{\bL^2}=\delta_{kj}.\label{Or}
\ee
From this, one knows that the following Poincare inequality holds:
\be
\sqrt{\l_1}\cdot\|\u\|_{\bL^2}\leq \|A^{\frac{1}{2}}\u\|_{\bL^2},\quad \forall\u\in\bH^1.\label{Poi}
\ee
Moreover, by (\ref{PI2}) and (\ref{Poi}) we have
\be
\|\u\|^2_{\infty}\leq C_0\cdot\|A\u\|_{\bL^2}\cdot\|\nabla\u\|_{\bL^2}.\label{Le4}
\ee
We have:
\bl\label{Le0} For $\epsilon>0$,
let $\cB_\epsilon:=\{\v\in\bH^1: \|\v\|_{\bH^1}\leq\epsilon\}$. Then $\cB_\epsilon$
is an absorbing set of $\{S(t); t\geq 0\}$, i.e., for any bounded set
$\cU\subset\bH^1$, there exists $t_\cU>0$ such that for any $t>t_\cU$
$$
S(t)\cU\subset\cB_\epsilon.
$$
\el
\begin{proof}
By the chain rule and (\ref{Poi}), we have
$$
\dif \|\u\|^2_{\bL^2}/\dif t=-2\nu\|\nabla\u\|^2_{\bL^2}
-2g^{\nu,\kappa}_N(\|\u\|^2_{\infty})\|\u\|^2_{\bL^2}\leq -2\nu\l_1\|\u\|^2_{\bL^2},
$$
which implies
\be
\|\u(t)\|^2_{\bL^2}\leq \|\u_0\|^2_{\bL^2}e^{-2\nu\l_1  t}.\label{PP5}
\ee
As the calculation of (\ref{Op11}), by  Young's inequality  we  have
\ce
\frac{\dif \|A^{\frac{1}{2}}\u\|^2_{\bL^2}}{\dif t}
&\leq&-\nu\|A\u\|^2_{\bL^2} +\frac{2\kappa N}{\nu}\|A^{\frac{1}{2}}\u\|^2_{\bL^2}\\
&\stackrel{(\ref{Le4})}{\leq}&-\nu\|A\u\|^2_{\bL^2}
+C_{\nu,\kappa} \cdot N\cdot\|A\u\|_{\bL^2}\cdot\|\u\|_{\bL^2}\\
&\leq&-\frac{\nu}{2}\|A\u\|^2_{\bL^2} +C_{\nu,\kappa} \cdot N^2\cdot\|\u\|^2_{\bL^2}\\
&\stackrel{(\ref{Poi})}{\leq}&-\frac{\nu\l_1}{2}\cdot\|A^{\frac{1}{2}}\u\|^2_{\bL^2}
+C_{\nu,\kappa} \cdot N^2\cdot\|\u_0\|^2_{\bL^2}\cdot e^{-2\nu\l_1  t}.
\de
Integrating this differential inequality yields that
\be
\|A^{\frac{1}{2}}\u(t)\|^2_{\bL^2}\leq e^{-\nu\l_1  t/2}\Big[\|A^{\frac{1}{2}}\u_0\|^2_{\bL^2}
+C_{\nu,\kappa}\cdot N^2\cdot\|\u_0\|^2_{\bL^2}\cdot(1-e^{-3\nu\l_1t/2})/(\nu\l_1)\Big].\label{PL1}
\ee
Hence, for any $\u_0\in\bH^1$
$$
\lim_{t\rightarrow\infty}\|S(t)\u_0\|^2_{\bH^1}=\lim_{t\rightarrow\infty}\|\u(t)\|^2_{\bH^1}=0.
$$
The result follows.
\end{proof}
We now use Lemma \ref{Le14} to prove the following compactness result.
\bl\label{Le00}
For any $t>0$, $S(t)$ is a compact operator from $\bH^1$ to $\bH^1$, i.e., maps a bounded set in $\bH^1$
into a relatively compact in $\bH^1$.
\el
\begin{proof}
Let $\cU\subset\bH^1$ be a bounded set.
Let $\Pi_n$ be the projection operator
from $\bH^1$ to span$\{\e_k: k=1,\cdots,n\}$, i.e.,
\be
\Pi_n\v:=\sum_{k=1}^n\<\v,\e_k\>_{\bL^2}\e_k.\label{Or1}
\ee

First of all, by (\ref{PL1}) we have
\be
\sup_{n\in\mN}\sup_{\u_0\in\cU}\|\Pi_nS(t)\u_0\|_{\bH^1}
\leq\sup_{\u_0\in\cU}\|S(t)\u_0\|_{\bH^1}<+\infty.\label{Po7}
\ee

Write
$$
\Pi^c_n:=I-\Pi_n
$$
By (\ref{Or}) and (\ref{Or1}) we have
$$
\Pi^c_n A=A\Pi^c_n.
$$
Thus, from (\ref{Op1}) we get
$$
\Pi^c_n\u(t)=\Pi^c_n\u_0-\nu\int^t_0A\Pi^c_n\u\dif s+\int^t_0\Pi^c_nB(\u,\u)\dif s
-\int^t_0g^{\nu,\kappa}_N(\|\u\|^2_\infty)\Pi^c_n\u\dif s.
$$

By the chain rule (cf. \cite[p.176, Lemma 1.2]{Te}) we have
\ce
\frac{\dif}{\dif t}\|\Pi^c_n\u(t)\|^2_{\bH^1}&=&-2\nu\|A\Pi^c_n\u\|^2_{\bL^2}
+2\<\Pi^c_nB(\u,\u),\Pi^c_n\u\>_{\bH^1}\\
&&-2g^{\nu,\kappa}_N(\|\u\|^2_\infty)\cdot\|\Pi^c_n\u\|^2_{\bH^1}\\
&\leq&-\nu\|A\Pi^c_n\u\|^2_{\bL^2}+\frac{1}{\nu}\|\Pi^c_nB(\u,\u)\|^2_{\bL^2}.
\de
Noting that
$$
\|\Pi^c_n\u\|^2_{\bH^1}=\|A^{\frac{1}{2}}\Pi^c_n\u\|^2_{\bL^2}\leq
\frac{1}{\lambda_n}\|A\Pi^c_n\u\|^2_{\bL^2}
$$
and
$$
\|\Pi^c_nB(\u,\u)\|^2_{\bL^2}\leq\|\u\|^2_{\infty}\cdot\|\nabla\u\|^2_{\bL^2}
\leq C\cdot\|A\u\|_{\bL^2}\cdot\|\nabla\u\|^3_{\bL^2}=:h(t),
$$
we have
\ce
\frac{\dif}{\dif t}\|\Pi^c_n\u(t)\|^2_{\bH^1}+\nu\lambda_n\|\Pi^c_n\u\|^2_{\bH^1}\leq
\frac{h(t)}{\nu}.
\de
Solving this differential inequality yields that
\ce
\|\Pi^c_n\u(t)\|^2_{\bH^1}&\leq& e^{-\nu\lambda_n t}\left(\|\Pi^c_n\u_0\|^2_{\bH^1}
+\int^t_0\frac{1}{\nu}e^{\nu\lambda_n s}h(s)\dif s\right)\\
&\leq& e^{-\nu\lambda_n t}\left(\|\u_0\|^2_{\bH^1}
+\frac{1}{\nu}\left(\int^t_0e^{2\nu\lambda_n s}\dif s\right)^{1/2}\left(\int^t_0h(s)^2\dif s\right)^{1/2}\right)\\
&\leq& e^{-\nu\lambda_n t}\|\u_0\|^2_{\bH^1}
+\frac{1}{\nu\sqrt{2\nu\lambda_n}}\left(\int^t_0h(s)^2\dif s\right)^{1/2}.
\de
On the other hand, by (\ref{Es222}) we have
$$
\int^t_0h(s)^2\dif s\leq \sup_{s\in[0,t]}\|\nabla\u(s)\|^6_{\bL^2}\int^t_0\|A\u\|^2_{\bL^2}\dif s
\leq \frac{1}{\nu}\Big(\frac{\kappa N}{\nu^2}\|\u_0\|^2_{\bL^2}+\|\nabla\u_0\|^2_{\bL^2}\Big)^4.
$$
Hence, by  $\lambda_n\uparrow\infty$ we obtain
$$
\lim_{n\rightarrow\infty}\sup_{\u_0\in\cU}\|\Pi^c_nS(t)\u_0\|^2_{\bH^1}
=\lim_{n\rightarrow\infty}\sup_{\u_0\in\cU}\|\Pi^c_n\u(t)\|^2_{\bH^1}=0,
$$
which combined with (\ref{Po7}) yields by  Lemma \ref{Le14} that
 $S(t)\cU$ is relatively compact in $\bH^1$.
\end{proof}

{\bf Proof of Theorem \ref{Th3}}:
It follows from \cite[p. 23 Theorem 1.1 and (1.12')]{Te2} and Lemmas \ref{Le0}
and \ref{Le00}.

\vspace{5mm}

{\bf Acknowledgements:}

The author would like to thank Professor Benjamin Goldys for
providing him an excellent environment to work in the University of New South Wales.
His work is supported by ARC Discovery grant DP0663153 of Australia.

\end{document}